\documentclass{article}
\usepackage[utf8]{inputenc}
\usepackage{graphicx}
\usepackage{times}
\usepackage{xcolor}
\usepackage{amsmath}
\usepackage{subfigure}
\usepackage{amssymb}
\usepackage{epstopdf}
\usepackage[numbers,sort&compress]{natbib}
\usepackage{appendix}
\usepackage{multirow}
\newcommand{\be}{\begin{equation}
\newcommand{\ee}{\end{equation}}}
\newcommand{\bea}{\begin{eqnarray}}
\newcommand{\eea}{\end{eqnarray}}
\newcommand{\nn}{\nonumber}

\begin{document}

\title{Solving conformable Gegenbauer differential equation and exploring its generating function  }

\author{ Mohamed Ghaleb Al-Masaeed$^{1}$, Eqab.M.Rabei$^{2}$, Sami I. Muslih$^{3, 4}$\\, and Dumitru Baleanu$^{5,6}$\\\\
$^1$Ministry of Education, Jordan\\
$^2$Physics Department, Faculty of Science, Al al-Bayt University,\\ P.O. Box 130040, Mafraq 25113, Jordan\\
$^3$Al-Azhar University-Gaza\\
$^4$Southern Illinois University,
Carbondale, IL, 62901 ~USA\\
$^5$ Department of Computer Science and Mathematics.\\
Lebanese American University, Beirut,Lebanon\\
$^6$ Institute of Space Sciences, Magurele–Bucharest, Romania\\\\ Email:moh.almssaeed@gmail.com\\eqabrabei@gmail.com\\sami.muslih@siu.edu\\dumitru.baleanu@lau.edu.lb}

\maketitle

%\begin{history}
%\received{Day Month Year}
%\revised{Day Month Year}
%\end{history}

\begin{abstract}
 In this manuscript, we address the resolution of conformable Gegenbauer differential equations. We demonstrate that our solution aligns precisely with the results obtained through the power series approach. Furthermore, we delve into the investigation and validation of various properties and recursive relationships associated with Gegenbauer functions. Additionally, we introduce and substantiate the conformable Rodriguez's formula and generating function\\

\textit{Keywords:}  conformable derivative, Gegenbauer differential equation, ultraspherical polynomials, 
\end{abstract}

\section{Introduction}
Special orthogonal polynomials started to appear in mathematics before the relevance of such a concept was made evident. Thus, whereas Legendre and Laplace used Legendre polynomials in celestial mechanics, Laplace used Hermite polynomials in his studies of probability \cite{mcbride1999special}. In 1784, Legendre made the first discovery of orthogonal polynomials \cite{legendre1785recherches}. Orthogonal polynomials are widely used to resolve ordinary differential equations when certain model constraints are met. Furthermore, orthogonal polynomials perform a significant role in the approximation theory \cite{bateman1953higher}. The class of classical orthogonal polynomials includes the polynomials Laguerre, Hermite, Legendre, and Gegenbauer. The Gegenbauer is one of the most significant classes of orthogonal polynomials. Since these polynomials are the foundation for many mathematical applications, scholars are interested in researching them in depth \cite{kazimouglu2023certain}.

 Gegenbauer polynomials belong to a distinct category of orthogonal polynomials. As elucidated in reference of \cite{kiepiela2003gegenbauer}, when employing conventional algebraic formulations,  a symbolic connection emerges between the generating function of Gegenbauer polynomials and the integral representation of real functions, typically denoted as  $T_R$.\\
 There are numerous applications of Gegenbauer polynomials that were used to solve the angular part of a second-order linear differential equation in dimensional space \cite{stillinger1977axiomatic,sadallah2006equations,palmer2004equations}.\\
Non-integer order calculus, sometimes known as fractional calculus, is a subject that extends from the Leibniz times and has recently received a lot of interest from mathematicians, physicists, and engineers \cite{agarwal2010concept,kai2010analysis,rabei2009fractional,herrmann2011fractional,rabei2007hamilton,kilbas2006theory,oldhamFractionalCalculusTheory1974,millerIntroductionFractionalIntegrals1993,podlubny1998fractional,uchaikin2013fractional}.  History has seen the introduction of several different forms of fractional derivative definitions, including Riemann-Louville, Caputo, Grunwald-Latnikov, Riesz, and Weyl. They all inherit nonlocal features from integrals because the majority of them are defined by fractional integrals. These concepts \cite{uchaikin2013fractional}, which differ from traditional Newton-Leibniz calculus and are significant in many application fields, frequently have the characteristics of heredity and nonlocality. The product rule, quotient rule, and chain rule for derivative operations are not followed by these derivatives. These contradictions cause problems and hassles in handling mathematics. Some academics developed the idea of local fractional derivative (LFD) to get around these problems \cite{zhao2017general}.We recommend that interested readers review these formulations that emerged by the end of the nineties to see \cite{kolwankar2013local}.

 In 2014, Khalil et al \cite{khalil2014new}, presented a new concept of derivative, called the conformable derivative.
 \\ \textbf{Definition 2.1.} Given a function $f\in [0,\infty) \to {R}$. The conformable  derivative of $f$ with order $\beta$ is defined by \cite{khalil2014new}

\be
\label{conformable}
T_\beta(f)(t)=\lim_{\epsilon \to 0}\frac{f(t+\epsilon t^{1-\beta})-f(t)}{\epsilon}
\ee
for all $t>0$, $\beta\in (0.1)$

 In this paper, we adopt $D^\beta f$ to denote the conformable  derivative (CD) of $f$ of order $\beta$, $ T_\beta(f)(t)$.\\
 The conformable derivative is considered a local fractional derivative \cite{teodoro2019review}. The properties of the traditional derivative are satisfied by this description \cite{abdeljawad2015conformable,atangana2015new}.  Khalil et al. answered the important question "What is the geometrical meaning of the conformable derivative?" This question was answered using the concept of fractional cords or conformable cords \cite{khalil2019geometric}. 
 
 The conformable calculus has several applications in special functions such as Legendre polynomials \cite{khader2017conformable,hammad2014legendre,
abul2020conformable}, associated Legendre polynomials \cite{shihab2021associated}, Laguerre polynomials \cite{abu2019laguerre,shat2019fractional},  associated  Laguerre polynomials \cite{doi:10.1142/S1793557123501814}, Bessel polynomials \cite{martinez2022novel,gökdoğan2015conformable}, hyper-geometric \cite{aldarawi2018conformable}, Chebyshev polynomials \cite{khader2017conformable,rababah2021conformable} and Hermite polynomials \cite{unal_solutions_2015}. In addition, these conformable special functions are used in many applications such as the associated Legendre polynomials \cite{shihab2021associated} and associated  Laguerre polynomials \cite{doi:10.1142/S1793557123501814} used to solve conformable Schrodinger equation for Coulomb potential \cite{al2022analytic}. Also, the Hermite polynomials \cite{unal_solutions_2015} used in the quantization of fractional harmonic oscillator using creation and annihilation operators \cite{al2021quantization}.
%%%%%%%%%%%%%%%%%%%%%%%%%%%%%%%%%%%%%%%%%%%%%%%%%%%%%%%%%%%%%%%%%%%%%%%%%%%%%%%%%%%%%%%%%%%%%%%%%%%%%%%%%%%%%%%%%%%%%%%%%%%%%%%%%%%%%%%%%%%%%%%%%%

%%%%%%%%%%%%%%%%%%%%%%%%%%%%%%%%%%%%%%%%%%%%%%%%%%%%%%%%%%%%%%%%%%%%%%%%%%%%%%
\section{Gegenbauer differential equation}
The  Gegenbauer differential equation takes the form \cite{doman2015classical}
\be
(1-x^{2}) y''-(2\lambda+1) x  y' +  n(n+2\lambda) y=0
\ee

The Gegenbauer polynomial  $C_n^{(\lambda)}$ of degree $n$ is defined by \cite{doman2015classical}
 \be
 C_n^{(\lambda)} = \sum_{k=0}^{[\frac{n}{2}]} (-1)^k \frac{\Gamma(n-k+\lambda)}{\Gamma(\lambda) k! (n-2k)!} (2z)^{n-2k}.
 \ee
 also the Gegenbauer polynomial is called Ultraspherical  polynomial \cite{elliott1960expansion}.\\
 The expression of $C_n^{(\lambda)}$ as a series \cite{doman2015classical}
\be
C_n^{(\lambda)} = (-1)^n \frac{\Gamma(\lambda+\frac{1}{2})\Gamma(n+2\lambda)}{2^n n!\Gamma(2\lambda) \Gamma(\lambda+n+\frac{1}{2})} (1-x^2)^{\frac{1}{2}-\lambda} \frac{d^n}{dx^n}[(1-x^2)^{n+\lambda-\frac{1}{2}}].
\ee
%%%%%%%%%%%%%%%%%%%%%%%%%%%%%%%%%%%%%%%%%%%%%%%%%%%%%%%%%%%%%%%%%%%%%%%%%%%%%%%%%%%%%%%%%%%%%%%%%%%%%%%%%%%%%%%%%%%%%%%%%%%%%%%%%%%%%%%%%%%
\section{  Gegenbauer Conformable Differential Equation }
We propose  the Gegenbauer conformable differential equation in the following way
\be
\label{CGEG}
(1-x^{2\alpha})D_x^\alpha D_x^\alpha y-\alpha(2\lambda+1) x^\alpha D_x^\alpha y + \alpha^2 n(n+2\lambda) y=0
\ee
when,\\
$\lambda=\frac{1}{2}$ this equation becomes conformable Legendre differential equation \cite{hammad2014legendre}.\\
$\lambda=0$ this equation becomes  Conformable Chebyshev differential equation of first kind \cite{rababah2021conformable}.\\
$\lambda=1$ this equation becomes  Conformable Chebyshev differential equation of the second kind. 
Let $y=\sum_{k=0}^\infty a_k x^{\alpha k} \to D_x^\alpha y=\sum_{k=1}^\infty \alpha k a_k x^{\alpha k-\alpha} \to  D_x^\alpha D_x^\alpha y=\sum_{k=2}^\infty \alpha^2 k(k-1) a_k x^{\alpha k-2\alpha} $, after substituting in eq.\eqref{CGEG}
\bea
\nn
\sum_{k=2}^\infty \alpha^2 k(k-1) a_k x^{\alpha k-2\alpha}&-& \sum_{k=2}^\infty \alpha^2 k(k-1) a_k x^{\alpha k}-\alpha^2(2\lambda+1)\sum_{k=1}^\infty  k a_k x^{\alpha k} \\&+&\alpha^2 \beta \sum_{k=0}^\infty a_k x^{\alpha k} =0
\eea
so,
\bea
\nn
 2\alpha^2 a_2&+& 6\alpha^2 a_3 x^\alpha+\sum_{k=2}^\infty \alpha^2 (k+1)(k+2) a_{k+2} x^{\alpha k}- \sum_{k=2}^\infty \alpha^2 k(k-1) a_k x^{\alpha k}\\\nn&-&\alpha^2(2\lambda+1)k a_1 x^\alpha-\alpha^2(2\lambda+1)\sum_{k=2}^\infty  k a_k x^{\alpha k}\\\nn &+& \alpha^2 \beta a_0+\alpha^2 \beta a_1 x^\alpha +\alpha^2 \beta \sum_{k=2}^\infty a_k x^{\alpha k} =0.
\eea
Thus, we have 
\be
a_2= -\frac{\beta}{2}a_0
\ee
\be
a_3= \frac{(2\lambda+1)-\beta}{6}a_1
\ee
\bea
\nn
a_{k+2} &=& \frac{k(k-1) +(2\lambda+1)k-\beta}{(k+1)(k+2)} a_k ,
\\&=&\frac{k(k+2\lambda)-\beta}{(k+1)(k+2)} a_k.
\eea
Let $k\to s-2 $, we get 
\bea
a_{s} &=&\frac{(s-2)(s-2+2\lambda)-\beta}{s(s-1)} a_{s-2}.
\eea
Let  $\beta = n(n+2\lambda)$ then we have 
\bea
\nn
a_{s} &=&\frac{(s-2)(s-2+2\lambda)-n(n+2\lambda)}{s(s-1)} a_{s-2},\\&=&\frac{(s-n-2)(s+n-2+2\lambda)}{s(s-1)} a_{s-2}
\eea
where $s\geq 2$
for $s=2$, we have 
\bea
\nn
a_{2} &=&\frac{-n(n+2\lambda)}{1\cdot2} a_{0}
\eea
for $s=3$, we have 
\bea
\nn
a_{3} &=&\frac{(1-n)(n+2\lambda+1)}{2\cdot3} a_{1}
\eea
for $s=4$, we have 
\bea
\nn
a_{4} &=&\frac{(2-n)(n+2\lambda+2)}{3\cdot4} a_{2}
\\\nn&=&\frac{n(n-2)(n+2\lambda)(n+2\lambda+2)}{1\cdot2\cdot3\cdot4} a_{0}
\eea
for $s=5$, we have 
\bea
\nn
a_{5} &=& \frac{(3-n)(n+2\lambda+3)}{4\cdot5} a_{3}\\\nn&=&
\frac{(n-1)(n-3)(n+2\lambda+1)(n+2\lambda+3)}{2\cdot3\cdot4\cdot5} a_{1}
\eea
for $s=6$, we have
\bea
\nn
a_{6} &=&\frac{(4-n)(n+2\lambda+4)}{5\cdot6} a_{4}
\\\nn&=&\frac{-n(n-2)(n-4)(n+2\lambda)(n+2\lambda+2)(n+2\lambda+4)}{1\cdot2\cdot3\cdot4\cdot5\cdot6} a_{0}
\eea
for $s=7$, we have
\bea
\nn
a_{7} &=& \frac{(5-n)(s+2\lambda+5)}{6\cdot7} a_{5}\\\nn&=&
\frac{-(n-1)(n-3)(n-5)(s+2\lambda+1)(s+2\lambda+3)(s+2\lambda+5)}{2\cdot3\cdot4\cdot5\cdot6\cdot7} a_{1}
\eea
The coefficient of $x^{\alpha s}$, $a_s=2^s\frac{\Gamma(\lambda+s)}{\Gamma(\lambda)s!}$to obtain the expression for $C_s^{(\lambda)}$as a polynomial \cite{doman2015classical}.
\be
\label{ser}
C_{\alpha s}^{(\lambda)}=\frac{1}{\Gamma(\lambda)} \sum_{s=0}^{\frac{n}{2}} \frac{(-1)^s \Gamma(\lambda+n-s)}{s! (n-2s)!} (2x^\alpha)^{n-2s}
\ee
where $s$ is the largest integer $\leq \frac{n}{2}$

%%%%%%%%%%%%%%%%%%%%%%%%%%%%%%%%%%%%%%%%%%%%%%%%%%%%%%%%%%%%%%%%%%%%%%%%%%%%%%
\subsection{Generating function }
\be
\label{gene1}
\sum_{n=0}^\infty t^{n\alpha} C_{\alpha n}^{(\lambda)}(x) = \frac{1}{(1-2x^\alpha t^\alpha+t^{2\alpha})^\lambda}
\ee
\textbf{Proof.} let 
\be
W(x^\alpha,t^\alpha)=\frac{1}{(1-2x^\alpha t^\alpha+t^{2\alpha})^\lambda}
\ee
According to reference \cite{doman2015classical} the following conformable differential equation can be considered
\bea
\label{w}
(1-x^{2\alpha})D_x^\alpha D_x^\alpha W-\alpha(2\lambda+1) x^\alpha D_x^\alpha W =- t^{\alpha(1-2\lambda)}D_t^\alpha [t^{\alpha(1+2\lambda)}D_t^\alpha W ].
\eea
Let us define 
\be
\sum_{n=0}^\infty t^{n\alpha} \Phi_n(x^\alpha) = \frac{1}{(1-2x^\alpha t^\alpha+t^{2\alpha})^\lambda},
\ee
after substituting left side in eq.\eqref{w}, we have 
\bea
\nn
&& \sum_{n=0}^\infty t^{n\alpha}[(1-x^{2\alpha})D_x^\alpha D_x^\alpha \Phi_n(x^\alpha)-\alpha(2\lambda+1) x^\alpha D_x^\alpha \Phi_n(x^\alpha)] \\&=&-\sum_{n=0}^\infty \Phi_n(x^\alpha) t^{\alpha(1-2\lambda)}D_t^\alpha [t^{\alpha(1+2\lambda)}D_t^\alpha t^{n\alpha} ].
\eea
so, we have 
\bea
\nn
&& \sum_{n=0}^\infty t^{n\alpha}[(1-x^{2\alpha})D_x^\alpha D_x^\alpha\Phi_n(x^\alpha)-\alpha(2\lambda+1) x^\alpha D_x^\alpha \Phi_n(x^\alpha) +  \alpha^2 n(n+2\lambda) \Phi_n(x^\alpha)]=0.
\eea
So, $\Phi(x^\alpha)$ satisfies conformable Gegenbauer’s equation. therefore 
$\Phi_n(x^\alpha)=C_{\alpha n}^{(\lambda)(x)}$. Thus, 
\be
\sum_{n=0}^\infty t^{n\alpha} C_{\alpha n}^{(\lambda)} (x)= \frac{1}{(1-2x^\alpha t^\alpha+t^{2\alpha})^\lambda}
\ee
when $\lambda=\frac{1}{2}$ this Generating function becomes conformable Legendre Generating function \cite{hammad2014legendre}.\\
%%%%%%%%%%%%%%%%%%%%%%%%%%%%%%%%%%%%%%%%%%%%%%%%%%%%%%%%%%%%%%%%%
\subsection{Differential relation}
One may  show that 
\be
\label{diff form}
D^{m\alpha}_x C_{\alpha n}^{(\lambda)} (x)= 2^m \alpha^m \lambda_m C_{\alpha (n-m)}^{(\lambda+m)} (x)
\ee
where $\lambda_m$ is the Pochammer symbol $\lambda_m=\lambda(\lambda+1)\cdots (\lambda+m-1)$.\\
\textbf{Proof}. If we differentiate the generating function eq.\eqref{gene1}, we have 
\bea
\sum_{n=0}^\infty t^{n\alpha} D^{\alpha}_x C_{\alpha s}^{(\lambda)}(x) = \frac{2\alpha \lambda t^\alpha}{(1-2x^\alpha t^\alpha+t^{2\alpha})^{\lambda+1}}= 2\alpha \lambda t^\alpha\sum_{n=0}^\infty t^{n\alpha} C_{\alpha n}^{(\lambda+1)}(x),
\eea
after equating powers of t on both sides, we obtain
\bea
 D^{\alpha}_x C_{\alpha s}^{(\lambda)}(x)= 2\alpha \lambda C_{\alpha (n-1)}^{(\lambda+1)}(x).
\eea
So, if differentiate the generating function $m$ times, we have eq. \eqref{diff form}.
Thus, if $\lambda=0$, eq.\eqref{diff form} becomes 
\be
\label{diff form000}
D^{m\alpha}_x C_{\alpha n}^{(0)}(x) = 2^m \alpha^m (m-1)! C_{\alpha (n-m)}^{(m)}(x), 
\ee
where $(\lambda)_{m}$ is Pochhammer Symbol and it is equal $(\lambda)_{m}=\frac{\Gamma(\lambda+m)}{\Gamma(\lambda)}$. 
%%%%%%%%%%%%%%%%%%%%%%%%%%%%%%%%%%%%%%%%%%%%%%%%%%%%%%%%%%%%%%%%%
\subsection{Recurrence relations }
Numerous recurrence relations can be derived using the generating function.\\
\textbf{The first relation }
 \be
 \label{drec1}
 (n+1)  C_{\alpha (n+1)}^{(\lambda)}(x)=2\lambda x^\alpha  C_{\alpha n}^{(\lambda+1)} (x)- 2\lambda  C_{\alpha (n-1)}^{(\lambda+1)}(x).
\ee
\textbf{Proof.} Differentiate the generating function in eq.\eqref{gene1}, we have 
\be
\label{rec1}
\sum_{n=1}^\infty n\alpha t^{n\alpha-\alpha} C_{\alpha n}^{(\lambda)}(x)=  \frac{2\alpha \lambda(x^\alpha - t^\alpha)}{(1-2x^\alpha t^\alpha+t^{2\alpha})^{\lambda+1}}.
\ee
Making use of equation \eqref{gene1}, we have 
\be
\sum_{n=1}^\infty n\alpha t^{n\alpha-\alpha} C_{\alpha n}^{(\lambda)}(x)=2\alpha \lambda x^\alpha \sum_{n=0}^\infty t^{n\alpha} C_{\alpha n}^{(\lambda+1)}(x) - 2\alpha \lambda \sum_{n=0}^\infty t^{n\alpha+\alpha} C_{\alpha n}^{(\lambda+1)}(x).
\ee
After equating the coefficient of $t^{n\alpha}$ , we obtain eq.\eqref{drec1}.\\
\textbf{The second relation }
\be
 \label{drec2}
(n+1)C_{\alpha (n+1)}^{(\lambda)}(x) = 2(n+\lambda)x^\alpha C_{\alpha n}^{(\lambda)}(x) - (n+2\lambda-1)C_{\alpha (n-1)}^{(\lambda)}(x)
\ee

\textbf{Proof.} Using eq.\eqref{rec1} and multiply it by $(1-2x^\alpha t^\alpha+t^{2\alpha})$, so, we have 
\be
(1-2x^\alpha t^\alpha+t^{2\alpha})\sum_{n=1}^\infty n\alpha t^{n\alpha-\alpha} C_{\alpha n}^{(\lambda)}(x)=  \frac{2\alpha \lambda(x^\alpha - t^\alpha)}{(1-2x^\alpha t^\alpha+t^{2\alpha})^{\lambda}}.
\ee
after substituting eq.\eqref{gene1}, we have 
\be
(1-2x^\alpha t^\alpha+t^{2\alpha})\sum_{n=1}^\infty n\alpha t^{n\alpha-\alpha} C_{\alpha n}^{(\lambda)}(x)=  2\alpha \lambda(x^\alpha - t^\alpha)\sum_{n=0}^\infty t^{n\alpha} C_{\alpha n}^{(\lambda)}(x).
\ee
Thus, we obtain 
\bea
(\sum_{n=1}^\infty n\alpha t^{n\alpha-\alpha} C_{\alpha n}^{(\lambda)}(x)&-&2x^\alpha \sum_{n=1}^\infty n\alpha t^{n\alpha} C_{\alpha n}^{(\lambda)}(x)+\sum_{n=1}^\infty n\alpha t^{n\alpha+\alpha} C_{\alpha n}^{(\lambda)}(x))
\\\nn&=&  2\alpha \lambda x^\alpha \sum_{n=0}^\infty t^{n\alpha} C_{\alpha n}^{(\lambda)}(x) - 2\alpha \lambda\sum_{n=0}^\infty t^{n\alpha+\alpha} C_{\alpha n}^{(\lambda)}(x).
\eea
After equating the coefficient of $t^{n\alpha}$ , we obtain  eq.\eqref{drec2}
%%%%%%%%%%%%%%%%%%%%%%%%%%%%%%%%%%%%%%%%%%%%%%%%%%%%%%%%%%%%%%%%%
\subsection{Conformable Rodrigues formula}
\be
\label{rodc}
C_{\alpha n}^{(\lambda)} = (-1)^n \frac{\Gamma(\lambda+\frac{1}{2})\Gamma(n+2\lambda)}{2^n \alpha^n n!\Gamma(2\lambda) \Gamma(\lambda+n+\frac{1}{2})} (1-x^{2\alpha})^{\frac{1}{2}-\lambda} D^{n\alpha}[(1-x^{2\alpha})^{n+\lambda-\frac{1}{2}}]
\ee
\textbf{Proof.} Let us define the conformable Gegenbauer polynomial
\be
\label{cr}
C_{\alpha n}^{(\lambda)}(x^\alpha)=A Q_n(x^\alpha)
\ee
Where $Q_n(x^\alpha)$ is $n th$ order conformable polynomial define by 
\be
\label{q}
Q_n(x^\alpha)=(1-x^{2\alpha})^{\frac{1}{2}-\lambda} D^{n\alpha}[(1-x^{2\alpha})^{n+\lambda-\frac{1}{2}}]
\ee
taking $ D^{n\alpha}[(1-x^{2\alpha})^{n+\lambda-\frac{1}{2}}]$ and rewrite it as 
\bea
\label{d}
 D^{n\alpha}[(1-x^{2\alpha})^{n+\lambda-\frac{1}{2}}]=(-1)^{n+\lambda-\frac{1}{2}} D^{n\alpha}[v^a u^a],
\eea
where $v=x^\alpha+1, u=x^\alpha-1, a=n+\lambda-\frac{1}{2}$ after using conformable Leibniz rule\cite{rabei_solution_2021}, we have 
\bea
D^{n\alpha}[v^a u^a]&=&\sum_{k=0}^n\left( \begin{array}{c}
      n\\
k      
\end{array}\right) D^{(n-k)\alpha} v^a D^{k\alpha} u^a,
\\\nn&=&\sum_{k=0}^n\left( \begin{array}{c}
      n\\
k      
\end{array}\right) \frac{a!\alpha^{n-k}}{(a-n+k)!} v^{a-n+k} \frac{a!\alpha^{k}}{(a-k)!} u^{a-k},
\\&=&\sum_{k=0}^n\left( \begin{array}{c}
      n\\
k      
\end{array}\right) \frac{(a!)^2\alpha^{n}}{(a-n+k)!(a-k)!} v^{a-n+k}  u^{a-k},
\eea
Substituting $ a=n+\lambda-\frac{1}{2}$, we obtain 
\bea
\nn
D^{n\alpha}[v^a u^a]&=&(vu)^{\lambda-\frac{1}{2}}\sum_{k=0}^n\left( \begin{array}{c}
      n\\
k      
\end{array}\right) \frac{((n+\lambda-\frac{1}{2})!)^2\alpha^{n}}{(\lambda+k-\frac{1}{2})!(n+\lambda-k-\frac{1}{2})!} v^{k}  u^{n-k},
\eea
substituting  in eq\eqref{d}, we have 
\bea
\nn
&& D^{n\alpha}[(1-x^{2\alpha})^{n+\lambda-\frac{1}{2}}]\\\nn&=&(-1)^{n} (1-x^{2\alpha})^{\lambda-\frac{1}{2}}\sum_{k=0}^n\left( \begin{array}{c}
      n\\
k      
\end{array}\right) \frac{((n+\lambda-\frac{1}{2})!)^2\alpha^{n}}{(\lambda+k-\frac{1}{2})!(n+\lambda-k-\frac{1}{2})!} (x^\alpha+1)^{k}  (x^\alpha-1)^{n-k},
\eea
Then, equation \eqref{q} will take the following form
\bea
Q_n(x^\alpha)= (-1)^{n} \sum_{k=0}^n\left( \begin{array}{c}
      n\\
k      
\end{array}\right) \frac{((n+\lambda-\frac{1}{2})!)^2\alpha^{n}}{(\lambda+k-\frac{1}{2})!(n+\lambda-k-\frac{1}{2})!} (x^\alpha+1)^{k}  (x^\alpha-1)^{n-k}.
\eea
equation  \eqref{cr} becomes 
\be
C_{\alpha n}^{(\lambda)}(x^\alpha)=A  (-1)^{n} \sum_{k=0}^n\left( \begin{array}{c}
      n\\
k      
\end{array}\right) \frac{((n+\lambda-\frac{1}{2})!)^2\alpha^{n}}{(\lambda+k-\frac{1}{2})!(n+\lambda-k-\frac{1}{2})!} (x^\alpha+1)^{k}  (x^\alpha-1)^{n-k}.
\ee
Taking  the limit of this equation when $k\to n$, we have 
\bea
\nn
C_{\alpha n}^{(\lambda)}(x^\alpha)&=&
A  (-1)^{n}  \frac{(n+\lambda-\frac{1}{2})!\alpha^{n}}{(\lambda-\frac{1}{2})!} (x^\alpha+1)^{n}
\\ &=& A  (-1)^{n}  \frac{\Gamma(n+\lambda+\frac{1}{2}) \alpha^{n}}{\Gamma(\lambda+\frac{1}{2})} (x^\alpha+1)^{n},
\eea
When $x^\alpha=1 \to C_{\alpha n}^{(\lambda)}(1)=\frac{\Gamma(2\lambda+n)}{\Gamma(2\lambda)n!}$. so,the constant A equal $A=\frac{\Gamma(2\lambda+n)\Gamma(\lambda+\frac{1}{2})}{(-2\alpha)^n\Gamma(2\lambda)n!\Gamma(n+\lambda+\frac{1}{2})}$. So, eq.\eqref{cr} becomes eq.\eqref{rodc}.
\bea
\nn
C_{\alpha 0}^{(\lambda)}(x^\alpha)&=&1.\\\nn
C_{\alpha 1}^{(\lambda)}(x^\alpha)&=&2 \lambda x^\alpha.\\\nn
C_{\alpha 2}^{(\lambda)}(x^\alpha)&=&\lambda_2  2 x^{2\alpha} -\lambda.\\\nn 
C_{\alpha 3}^{(\lambda)}(x^\alpha)&=&\lambda_3\frac{4}{3}x^{3\alpha}-\lambda_2 2x^\alpha.\\\nn
C_{\alpha 4}^{(\lambda)}(x^\alpha)&=& \lambda_4  \frac{2}{3} x^{4\alpha}- \lambda_3 2 x^{2\alpha}+\lambda_2. \\\nn
C_{\alpha 5}^{(\lambda)}(x^\alpha)&=&\lambda_5 \frac{4x^{5\alpha}}{15}-\lambda_4 \frac{4x^{3\alpha}}{3}+\lambda_3 x^\alpha,
\eea
where $\lambda_n$ is the Pochammer symbol $\lambda_n=\lambda(\lambda+1)\cdots (\lambda+n-1)$.
\begin{figure}[htb!]
    \centering
    \includegraphics[width=0.68\textwidth]{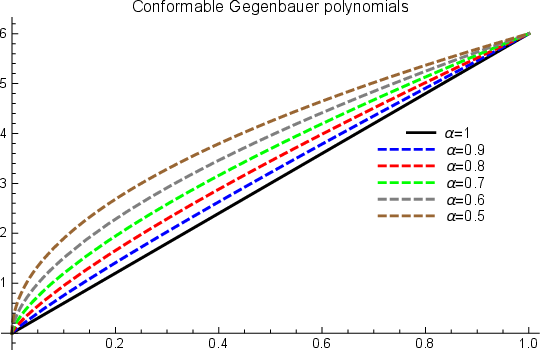}
    \caption{Plot of $C_{\alpha 1}^{(\lambda)}(x^\alpha)=2 \lambda x^\alpha.$ with different values of $\alpha$ when $\lambda=3$.}
    \label{fig:enter-label}
\end{figure}
\begin{figure}[htb!]
    \centering
    \includegraphics[width=0.68\textwidth]{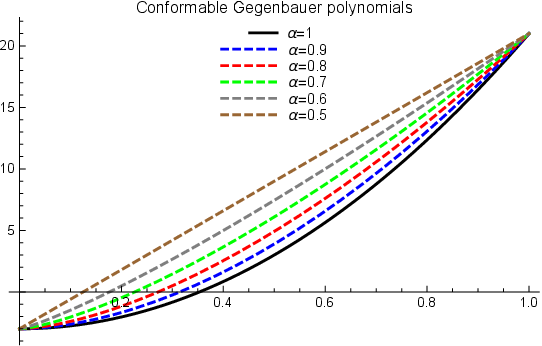}
    \caption{Plot of $C_{\alpha 2}^{(\lambda)}(x^\alpha)=\lambda_2  2 x^{2\alpha} -\lambda.$ with different values of $\alpha$ when $\lambda=3$.}
    \label{fig:enter-label}
\end{figure}
\begin{figure}[htb!]
    \centering
    \includegraphics[width=0.68\textwidth]{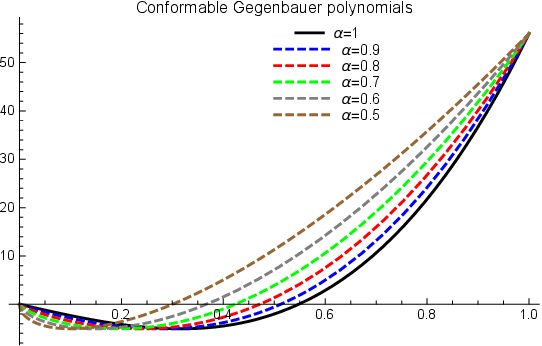}
    \caption{Plot of $C_{\alpha 3}^{(\lambda)}(x^\alpha)=\lambda_3\frac{4}{3}x^{3\alpha}-\lambda_2 2x^\alpha.$ with different values of $\alpha$ when $\lambda=3$.}
    \label{fig:enter-label}
\end{figure}
\clearpage
\begin{figure}[htb!]
    \centering
    \includegraphics[width=0.68\textwidth]{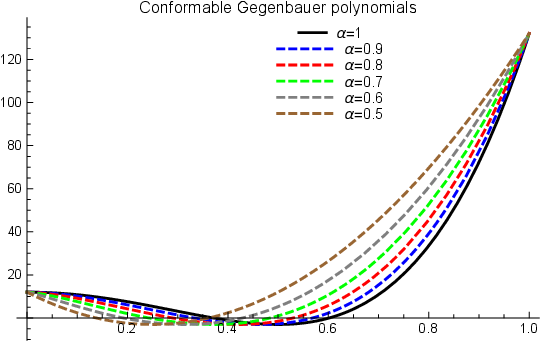}
    \caption{Plot of $C_{\alpha 4}^{(\lambda)}(x^\alpha)= \lambda_4  \frac{2}{3} x^{4\alpha}- \lambda_3 2 x^{2\alpha}+\lambda_2$ with different values of $\alpha$ when $\lambda=3$.}
    \label{fig:enter-label}
\end{figure}

\begin{figure}[htb!]
    \centering
    \includegraphics[width=0.68\textwidth]{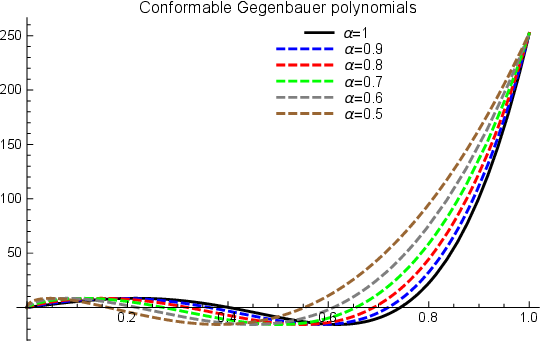}
    \caption{Plot of $C_{\alpha 5}^{(\lambda)}(x^\alpha)=\lambda_5 \frac{4x^{5\alpha}}{15}-\lambda_4 \frac{4x^{3\alpha}}{3}+\lambda_3 x^\alpha$ with different values of $\alpha$ when $\lambda=3$.}
    \label{fig:enter-label}
\end{figure}

when,
$\lambda=\frac{1}{2}$ the eq.\eqref{rodc} becomes conformable Legendre polynomial \cite{hammad2014legendre}.
\be
C_{\alpha n}^{(\frac{1}{2})} =  \frac{1}{2^n \alpha^n n! }  D^{n\alpha}[(x^{2\alpha}+1)^{n}]=P_{n\alpha}(x)
\ee
$\lambda=0$ the eq.\eqref{rodc} becomes  Conformable Chebyshev polynomial of first kind \cite{rababah2021conformable}.
\be
C_{\alpha n}^{(0)} = (-1)^n \frac{2^n n! }{ \alpha^n (2n)! } (1-x^{2\alpha})^{\frac{1}{2}} D^{n\alpha}[(1-x^{2\alpha})^{n-\frac{1}{2}}]=T_{n\alpha}(x)
\ee

%%%%%%%%%%%%%%%%%%%%%%%%%%%%
\subsection{Orthogonality and normalization}
The conformable  Gegenbauer polynomials $C_{\alpha n}^{(\lambda)}(x^\alpha)$ are satisfying   the orthogonality relation 
\be
\label{orth}
\int_{-1}^1  (1-x^{2\alpha})^{\lambda-\frac{1}{2}} C_{\alpha m}^{(\lambda)}C_{\alpha n}^{(\lambda)}(x^\alpha) d^\alpha x=0.
\ee
\textbf{Proof.} We can be written  the conformable Gegenbauer differential equation \eqref{CGEG} for $C_{\alpha n}^{(\lambda)}(x^\alpha)$ in  Conformable Sturm-Liouville form \cite{al-refai_fundamental_2017} 
\be
\label{cn}
D_x^\alpha[(1-x^{2\alpha})^{\lambda+\frac{1}{2}}D_x^\alpha C_{\alpha n}^{(\lambda)}(x^\alpha)] =-\alpha^2 n(n+2\lambda) (1-x^{2\alpha})^{\lambda-\frac{1}{2}}C_{\alpha n}^{(\lambda)}(x^\alpha).
\ee
Replace $n$ by $m $ and rewrite the above equation, we have 
\be
\label{cm}
D_x^\alpha[(1-x^{2\alpha})^{\lambda+\frac{1}{2}}D_x^\alpha C_{\alpha m}^{(\lambda)}(x^\alpha)] =-\alpha^2 m(m+2\lambda) (1-x^{2\alpha})^{\lambda-\frac{1}{2}} C_{\alpha m}^{(\lambda)}(x^\alpha).
\ee
After multiplying eq.\eqref{cn} by $C_{\alpha m}^{(\lambda)}(x^\alpha)$ and eq.\eqref{cm} by $C_{\alpha n}^{(\lambda)}(x^\alpha)$ and subtracting these equations, we obtain
\bea
\nn
&&C_{\alpha m}^{(\lambda)}(x^\alpha)D_x^\alpha[(1-x^{2\alpha})^{\lambda+\frac{1}{2}}D_x^\alpha C_{\alpha n}^{(\lambda)}(x^\alpha)]- C_{\alpha n}^{(\lambda)}(x^\alpha)D_x^\alpha[(1-x^{2\alpha})^{\lambda+\frac{1}{2}}D_x^\alpha C_{\alpha m}^{(\lambda)}(x^\alpha)]
\\&=&\alpha^2 [m(m+2\lambda)-n(n+2\lambda)] (1-x^{2\alpha})^{\lambda-\frac{1}{2}} C_{\alpha m}^{(\lambda)}C_{\alpha n}^{(\lambda)}(x^\alpha).
\eea
Integrating both sides of this equation, we have 
\bea
\nn
&&\int_{-1}^1 C_{\alpha m}^{(\lambda)}(x^\alpha)D_x^\alpha[(1-x^{2\alpha})^{\lambda+\frac{1}{2}}D_x^\alpha C_{\alpha n}^{(\lambda)}(x^\alpha)]d^\alpha x-\int_{-1}^1 C_{\alpha n}^{(\lambda)}(x^\alpha)D_x^\alpha[(1-x^{2\alpha})^{\lambda+\frac{1}{2}}D_x^\alpha C_{\alpha m}^{(\lambda)}(x^\alpha)]d^\alpha x
\\&=&\alpha^2 [m(m+2\lambda)-n(n+2\lambda)]\int_{-1}^1  (1-x^{2\alpha})^{\lambda-\frac{1}{2}} C_{\alpha m}^{(\lambda)}C_{\alpha n}^{(\lambda)}(x^\alpha) d^\alpha x.
\eea
using the integration by parts for the left side of this equation \cite{abdeljawad2015conformable}, we have 
\bea
\nn
&&\int_{-1}^1 (1-x^{2\alpha})^{\lambda+\frac{1}{2}} D_x^\alpha C_{\alpha m}^{(\lambda)}(x^\alpha)D_x^\alpha C_{\alpha n}^{(\lambda)}(x^\alpha)d^\alpha x-\int_{-1}^1 (1-x^{2\alpha})^{\lambda+\frac{1}{2}} D_x^\alpha C_{\alpha m}^{(\lambda)}(x^\alpha)D_x^\alpha C_{\alpha n}^{(\lambda)}(x^\alpha)d^\alpha x
\\&=&\alpha^2 [m(m+2\lambda)-n(n+2\lambda)]\int_{-1}^1  (1-x^{2\alpha})^{\lambda-\frac{1}{2}} C_{\alpha m}^{(\lambda)}C_{\alpha n}^{(\lambda)}(x^\alpha) d^\alpha x.
\eea
For $m \neq n$, the orthogonality relation \eqref{orth} is satisfied.\\
The normalization of conformable  Gegenbauer polynomials $C_{\alpha n}^{(\lambda)}(x^\alpha)$ are satisfies in case $m =n$
\be
\label{norm}
\int_{-1}^1  (1-x^{2\alpha})^{\lambda-\frac{1}{2}} [C_{\alpha n}^{(\lambda)}]^2  d^\alpha x= 2^{1-2\lambda}\alpha^{-\frac{2}{\alpha}} \frac{\Gamma(n+2\lambda)\Gamma(\lambda+n)\Gamma(\frac{5}{2}-\alpha-\frac{1}{\alpha})\Gamma(n+\lambda+\frac{3}{2}-\frac{1}{\alpha})}{ n! [\Gamma(\lambda)]^2\Gamma(\lambda+n+\frac{1}{2})\Gamma(n+\lambda+2-\alpha)}.
\ee

\textbf{Proof.} using the eq.\eqref{rodc}, we have 
\be
\label{general ss}
\int_{-1}^1  (1-x^{2\alpha})^{\lambda-\frac{1}{2}} [C_{\alpha n}^{(\lambda)}]^2  d^\alpha x=  \frac{(-1)^n\Gamma(\lambda+\frac{1}{2})\Gamma(n+2\lambda)}{2^n \alpha^n n!\Gamma(2\lambda) \Gamma(\lambda+n+\frac{1}{2})}\int_{-1}^1   C_{\alpha n}^{(\lambda)}D^{n\alpha}[(1-x^{2\alpha})^{n+\lambda-\frac{1}{2}}]  d^\alpha x,
\ee
using the integration by parts for the right side of this equation \cite{abdeljawad2015conformable}, we have 
\be
\label{cint}
\int_{-1}^1   C_{\alpha n}^{(\lambda)}D^{n\alpha}[(1-x^{2\alpha})^{n+\lambda-\frac{1}{2}}]  d^\alpha x=(-1)^n \int_{-1}^1   (1-x^{2\alpha})^{n+\lambda-\frac{1}{2}} D^{n\alpha}C_{\alpha n}^{(\lambda)}  d^\alpha x
\ee
using the eq.\eqref{diff form}, when $m =n$, we have 
\be
D^{n\alpha}_x C_{\alpha n}^{(\lambda)} (x)= 2^n \alpha^n \lambda_n C_{\alpha (0)}^{(\lambda+m)} (x)= 2^n \alpha^n \frac{\Gamma(\lambda+n)}{\Gamma(\lambda)}.
\ee
So, the integration \eqref{cint} becomes 
\be
\label{c2222}
\int_{-1}^1   C_{\alpha n}^{(\lambda)}D^{n\alpha}[(1-x^{2\alpha})^{n+\lambda-\frac{1}{2}}]  d^\alpha x=(-1)^n  2^n \alpha^n \frac{\Gamma(\lambda+n)}{\Gamma(\lambda)}\int_{-1}^1   (1-x^{2\alpha})^{n+\lambda-\frac{1}{2}}  d^\alpha x.
\ee
using conformable Gamma and Beta functions \cite{sarikaya2020some}, the solution of this integral is given by 
\bea
\int_{-1}^1   (1-x^{2\alpha})^{n+\lambda-\frac{1}{2}}  d^\alpha x&=&  \frac{\Gamma^\alpha(\frac{3}{2}\alpha-\alpha^2)\Gamma^\alpha(\alpha(n+\lambda+\frac{1}{2}))}{\alpha \Gamma^\alpha(\alpha[n+\lambda+1-\alpha]+1)},
\\\nn&=& \alpha^{\frac{1}{2}-\frac{2}{\alpha}}\frac{\Gamma(\frac{5}{2}-\alpha-\frac{1}{\alpha})\Gamma(n+\lambda+\frac{3}{2}-\frac{1}{\alpha})}{\Gamma(n+\lambda+2-\alpha)},
\eea
where the conformable gamma function is defined as  $\Gamma^\alpha(p)=\alpha^{\frac{p+\alpha-1}{\alpha}} \Gamma(\frac{p+\alpha-1}{\alpha})$\cite{sarikaya2020some} 
So, the integral \eqref{c2222} becomes 
\be
\int_{-1}^1   C_{\alpha n}^{(\lambda)}D^{n\alpha}[(1-x^{2\alpha})^{n+\lambda-\frac{1}{2}}]  d^\alpha x=(-1)^n 2^n \alpha^n \frac{\Gamma(\lambda+n)}{\Gamma(\lambda)}\alpha^{\frac{1}{2}-\frac{2}{\alpha}}\frac{\Gamma(\frac{5}{2}-\alpha-\frac{1}{\alpha})\Gamma(n+\lambda+\frac{3}{2}-\frac{1}{\alpha})}{\Gamma(n+\lambda+2-\alpha)}
\ee
After substituting in eq.\eqref{general ss}, we have 
\be
\nn
\int_{-1}^1  (1-x^{2\alpha})^{\lambda-\frac{1}{2}} [C_{\alpha n}^{(\lambda)}]^2  d^\alpha x= \alpha^{\frac{1}{2}-\frac{2}{\alpha}} \frac{\Gamma(\lambda+\frac{1}{2})\Gamma(n+2\lambda)}{ n!\Gamma(2\lambda) \Gamma(\lambda+n+\frac{1}{2})}\frac{\Gamma(\lambda+n)}{\Gamma(\lambda)}\frac{\Gamma(\frac{5}{2}-\alpha-\frac{1}{\alpha})\Gamma(n+\lambda+\frac{3}{2}-\frac{1}{\alpha})}{\Gamma(n+\lambda+2-\alpha)}.
\ee
Using the gamma function properties, we have eq.\eqref{norm}

%%%%%%%%%%%%%%%%%%%%%%%%%%%%%%%%%%%%%%%%%%%%%%%%%%%%%%%%%%%%%%%%%%%%%%%%%%%%%%%%%%%%%%%%%%%%%%%%%%%%%%
\section{Another formula for Conformable Gegenbauer equation  }
The Gegenbauer differential equation for an integer $n$ has solutions in the form of Gegenbauer polynomials $T_{\alpha n}^\beta(x)$. They are generalizations of the associated Legendre polynomials to $2(\beta+1)$-D space and have a proportion to the ultraspherical polynomials \cite{geg},  when $\lambda=\beta +\frac{1}{2}$, the eq.\eqref{CGEG} becomes 
\be
\label{CGEG2}
D_x^\alpha  D_x^\alpha y-2\alpha(\beta+1) x^\alpha D_x^\alpha y + \alpha^2 n(n+2\beta+1) y=0.
\ee
So, the series solution in eq.\eqref{ser} becomes 
\be
\label{ser2}
T_{\alpha n}^\beta(x)=\frac{1}{\Gamma(\beta + \frac{1}{2})} \sum_{s=0}^{\frac{n}{2}} \frac{(-1)^s \Gamma(\beta+n-s + \frac{1}{2})}{s! (s-2n)!} (2x^\alpha)^{s-2n}.
\ee
The solution with Rodrigues formula in eq.\eqref{rodc}becomes
\be
\label{rodc2}
T_{\alpha n}^\beta(x) = (-1)^n \frac{\Gamma(n+2\beta+1)}{2^{n+\beta} \alpha^n n! \Gamma(n+\beta+1)} (1-x^{2\alpha})^{-\beta} D^{n\alpha}[(1-x^{2\alpha})^{n+\beta}]
\ee
The ordinary form of these equations \eqref{CGEG2},\eqref{ser2}, and \eqref{rodc2} is recovered when  $\alpha=1$ in ref \cite{morse1954methods}
%%%%%%%%%%%%%%%%%%%%%%%%%%%%

%%%%%%%%%%%%%%%%%%%%%%%%%%%%%%%%%%%%%%%%%%%%%%%%%
\section{Conclusions}
We solved our proposed Gegenbauer conformable differential equation and we obtained the conformable Gegenbauer generating function. We observed this generating function goes to conformable generating function when  $\lambda=\frac{1}{2}$. In addition, we proved some recurrence relations for the conformable Gegenbauer function. And we have seen the conformable Gegenbauer function satisfies  the orthonormality relations. We painted the conformable Gegenbauer polynomials for different values of $n $ with different values of $\alpha$, and we found the curves gradually convert to $\alpha=1$ for each $n$.
%%%%%%%%%%%%%%%%%%%%%%%%%%%%%%%%%%%%%%%%%%%%%%%%%%%%%%%%%%%%%%%%%%%%%%%%%%%%%%%%%%%%%%%%%%%%%%%%%%%%%%%
\bibliography{ref} 
\bibliographystyle{unsrtnat}
\end{document}